\documentclass[12pt]{article}
\usepackage{amsmath,amssymb,amsbsy,amsfonts,amsthm,latexsym,
            amsopn,amstext,amsxtra,euscript,amscd,amsthm,url}

\usepackage[usenames, dvipsnames]{color}
\usepackage[utf8]{inputenc}
\usepackage[english]{babel}
\usepackage{mathrsfs}
\usepackage{caption}

\newtheorem{theorem}{Theorem}
\newtheorem{remark}{Remark}
\newtheorem{definition}{Definition}

\newtheorem{corol}{Corollary}

\newcommand{\icomp}{\texttt{i}}
\newcommand{\rd}{\, {\rm d}}

\DeclareMathOperator{\e}{e}
\DeclareMathOperator{\id}{id}
\DeclareMathOperator{\lcm}{lcm}

\begin{document}
\date{}

\title{On the multivariable generalization of Anderson-Apostol sums}
\author{Isao Kiuchi, Friedrich Pillichshammer, and Sumaia Saad Eddin}

\maketitle
{\def\thefootnote{}
\footnote{{\it Mathematics Subject Classification 2010: 11A25, 11N37, 11Y60.}}

\begin{abstract}
In this paper, we give various identities for the weighted average of the product of generalized Anderson-Apostol sums with weights concerning completely multiplicative function, completely additive function, logarithms, the Gamma function, Bernoulli polynomials and binomial coefficients respectively. Our results generalize many known results.
\end{abstract}
\maketitle
\section{Introduction}  

Let $\mathbb{N}$ be the set of positive integers and $\mathbb{N}_0=\mathbb{N} \cup\{0\}$. For $j,k \in \mathbb{N}$, let $\gcd (j, k)$ and $\lcm(j,k)$ denote their greatest common divisor and least common multiple, respectively. The Ramanujan sum $c_{k}$ for $k \in \mathbb{N}$ is an arithmetic function which is defined as
$$ 
c_{k}(j)=\sum_{d|\gcd(j,k)}d\mu\left(\frac{k}{d}\right)=\sum_{\substack{m=1\\ \gcd(m,k)=1}}^{k}\exp \left( 2\pi \icomp m j/k\right), \ \ \ \mbox{for $j \in \mathbb{N}$,} 
$$  
where $\mu$ is the M\"obius function and $\icomp=\sqrt{-1}$. This arithmetic function has received attention due to its connections to many problems in number theory, such as the proof of Vinogradov's theorem which implies that any sufficiently large odd integer can be written as a sum of three prime numbers.

Nowadays various generalizations of the Ramanujan sum are studied. One of these is according to Cohen \cite{C1,C2,C3} who defined the arithmetic function $c_k^{(a)}$ for $k,a \in \mathbb{N}$ by 
\begin{equation}\label{Cohen}
c_{k}^{(a)}(j) := \sum_{\substack{d |k\\d^a|j}}d^{a}\mu\left(\frac{k}{d}\right)=\sum_{\substack{m=1\\ (m,k^a)_a=1}}^{k^a}\exp \left( 2\pi \icomp m j/k^a\right),  \ \ \ \mbox{for $j \in \mathbb{N}$,} 
\end{equation} 
where $(m,k^a)_a=1$ means that no positive integer greater than one satisfies $d|k$ and $d^a|m$. 
Many interesting properties and useful formulas of Ramanujan-Cohen sum Eq.~\eqref{Cohen} have been given by McCarthy~\cite{Mc1,Mc3,Mc2}, K\"uhn and Robles~\cite{K.R1}, Robles and Roy~\cite{R.R1}, and others.  

In 2006, Mednykh and Nedela~\cite{M.N1} studied the average of products of Ramanujan sums with the aim to handle certain problems of enumerative combinatorics. To this end they established the function $E$ of $n$ variables given by
\begin{equation}\label{MN1}
E(\boldsymbol{k})=\frac{1}{[\boldsymbol{k}]}\sum_{j=1}^{[\boldsymbol{k}]}c_{k_1}(j)\cdots c_{k_n}(j), \ \ \ \mbox{for $\boldsymbol{k} \in \mathbb{N}^n$,} 
\end{equation}
where here and throughout this paper  we write $\boldsymbol{k}=(k_1, k_2, \ldots, k_n)$ and $[\boldsymbol{k}]=\lcm(k_1, \ldots, k_n)$. 
    Furthermore, for later use we mention already here that if integers $d_1|k_1,\ldots,d_n|k_n$, then we will write $\boldsymbol{d}|\boldsymbol{k}$, $\boldsymbol{f}(\boldsymbol{d})=f_1(d_1)\cdots f_n(d_n),$ and $\boldsymbol{f}(d)=f_1(d)\cdots f_n(d),$ (similarly for $g$ and $h$).

It is shown in \cite{M.N1} that $E(\boldsymbol{k})$ belongs to $\mathbb{N}$ for all $\boldsymbol{k}\in \mathbb{N}^n$. Arithmetic properties of $E$ have been studied by Liskovets~\cite{L1}, who called $E$ the ``{\it orbicyclic}" arithmetic function, and by  T\'oth~\cite{To1}. In particular, T\'oth proved that $E$ is a multiplicative function which can also be represented in the form 
\begin{equation*}
E(\boldsymbol{k})=\sum_{\boldsymbol{d}|\boldsymbol{k}}\frac{d_1\mu(k_1/d_1)\cdots d_n\mu(k_n/d_n)}{[\boldsymbol{d}]}.
\end{equation*}

Recently, T\'oth~\cite{To2} extended the definition of $E$ in Eq.~\eqref{MN1} and considered the weighted average
\begin{equation}\label{toav1}
S_r(\boldsymbol{k}):=\frac{1}{[\boldsymbol{k}]^{r+1}}\sum_{j=1}^{[\boldsymbol{k}]}j^r c_{k_1}(j)\cdots c_{k_n}(j),
\end{equation}
where $r\in \mathbb{N}_0$. Clearly, $S_0=E$. He proved that
\begin{equation}
\label{resth1}
S_r(\boldsymbol{k})= \frac{1}{2[\boldsymbol{k}]}\prod_{j=1}^n \phi(k_j)+\frac{1}{r+1}\sum_{m=0}^{\lfloor r/2\rfloor}\binom{r+1}{2m}\, \frac{B_{2m}}{[\boldsymbol{k}]^{2m}}\, g_m(\boldsymbol{k}),   
\end{equation}
where $\phi$ and $B_m$ with $m\in \mathbb{N}_0$ are the Euler totient function and Bernoulli numbers, respectively, and the functions $g_m(\boldsymbol{k})$ are given by 
\begin{equation*}
g_m(\boldsymbol{k})=\sum_{\boldsymbol{d}|\boldsymbol{k}}\frac{d_1\mu(k_1/d_1)\cdots d_n\mu(k_n/d_n)}{[\boldsymbol{d}]^{1-2m}}.
\end{equation*}
Furthermore, for $r=1$, it is shown that 
\begin{equation*}
S_1(\boldsymbol{k})= \frac{1}{2[\boldsymbol{k}]}\prod_{j=1}^n \phi(k_j)+\frac{E(\boldsymbol{k})}{2}.
\end{equation*}

Another generalization of the Ramanujan sum $c_k(j)$ has been given by Anderson and Apostol \cite{AA} who introduced the arithmetic function $s_{k}$ for $k \in \mathbb{N}$ defined by the identity   
\begin{equation}\label{A.A}
 s_{k}(j)=\sum_{d|\gcd(k,j)}f(d)g\left(\frac{k}{d}\right), \ \ \ \mbox{for $j \in \mathbb{N}$,} 
\end{equation}
with arithmetic functions $f$ and $g$. Then one can generalize the function $E(\boldsymbol{k})$ in Eq.~\eqref{MN1} to
\begin{equation}
\widetilde{E}(\boldsymbol{k})=\frac{1}{[\boldsymbol{k}]}\sum_{j=1}^{[\boldsymbol{k}]}s_{k_1}(j)\cdots s_{k_n}(j).
\end{equation}

In~\cite{To1}, T\'oth investigated this function and proved that if $f$ and $g$ are multiplicative, then the function $\widetilde{E}$ is multiplicative as well. 
 
Generalizing Eq.~\eqref{toav1} we obtain a weighted average of products of $s_{k_i}$ of the form 
 \begin{equation}
\widetilde{S}_r(\boldsymbol{k}):=\frac{1}{[\boldsymbol{k}]^{r+1}}\sum_{j=1}^{[\boldsymbol{k}]}j^r \ s_{k_1}(j)\cdots s_{k_n}(j).
\end{equation}
Recently, Ikeda, Kiuchi and Matsuoka~\cite{I.K.M} proved that
$$
\widetilde{S}_r(\boldsymbol{k})= \frac{1}{2[\boldsymbol{k}]} \prod_{j=1}^n f_j \ast g_j(k_j) +\frac{1}{r+1}\sum_{m=0}^{\lfloor r/2\rfloor}\binom{r+1}{2m}\ \frac{B_{2m}}{[\boldsymbol{k}]^{2m}} \ \widetilde{g}_m(\boldsymbol{k}),   
$$
where
$$\widetilde{g}_m(\boldsymbol{k})=\sum_{\boldsymbol{d}|\boldsymbol{k}}\frac{\boldsymbol{f}(\boldsymbol{d})}{[\boldsymbol{d}]^{1-2m}}\ \boldsymbol{g}\left(\frac{\boldsymbol{k}}{\boldsymbol{d}}\right),$$
and the symbol $*$ denotes the Dirichlet convolution of arithmetic functions $$g \ast f (n)=\sum_{d|n} f(d) g\left(\frac{n}{d}\right),$$
Compare this result with Eq.~\eqref{resth1}.\\ 

We propose a further generalization of Anderson-Apostol sums which will be given in the following definition:

\begin{definition}
For fixed $a\in \mathbb{N}$ and for given arithmetic functions $f$, $g$ and $h$, define the functions $s^{(a)}_{f,g,h}$ by 
\begin{equation}
s^{(a)}_{f, g, h}(k, j)=\sum_{\substack{d|k\\d^a|j}}f(d)g\left(\frac{k}{d}\right)h\left(\frac{j}{d^a}\right) \ \ \ \ \mbox{ for $k,j \in \mathbb{N}$.}
\end{equation}
We call these functions {\it generalized Anderson-Apostol sums}.
\end{definition}

Let the arithmetic function ${\bf 1}: \mathbb{N} \rightarrow \mathbb{N}$ be given by ${\bf 1}(n)=1$ for all $n$. It is easily to see that $s^{(1)}_{f, g, \bf 1}(k, j)=s_k(j)$ and $s^{(a)}_{\id_a, \mu, \bf 1}(k, j)=c_k^{(a)}(j)$, where  $\id_r(j):=j^r$. In~\cite{K2017}, the first author investigated the weighted average of $s^{(a)}_{f, g, {\bf 1}}(k,j)=s^{(a)}_{k}(j)$ with the weight function $\id_r$. He proved that
 \begin{equation*}
     \frac{1}{k^{ar}}\sum_{j=1}^{k^a}j^rs^{(a)}_{k}(j)=
     \frac{1}{2}f*g (k)+\frac{1}{r+1}\sum_{m=0}^{\lfloor r/2\rfloor}\binom{r+1}{2m} B_{2m}\left(f*\id_{a(1-2m)}\cdot g \right)(k).
 \end{equation*}

Now it is the aim of the present paper to study weighted averages of products  $$s^{(a)}_{f_1,g_1, h_1}(k,j) \cdots s^{(a)}_{f_n, g_n, h_n}(k, j)$$ for finite sequences of given arithmetic functions $(f_j)_j$, $(g_j)_j$, and $(h_j)_j$, with weight functions $\omega$ that are either completely multiplicative or completely additive. Put
\begin{equation}\label{defU}
    U_{\omega}^{(a)}(\boldsymbol{k}):=\sum_{j=1}^{[\boldsymbol{k}]^a}\omega(j)s^{(a)}_{f_1, g_1, h_1}(k_1,j)\cdots s^{(a)}_{f_n, g_n,h_n}(k_n, j).
\end{equation}
If $h_1=\ldots=h_n={\bf 1}$, then we write 
\begin{equation*}
    \widetilde{U}_{\omega}^{(a)}(\boldsymbol{k}):=\sum_{j=1}^{[\boldsymbol{k}]^a}\omega(j)s^{(a)}_{f_1, g_1, {\bf 1}}(k_1,j)\cdots s^{(a)}_{f_n, g_n,{\bf 1}}(k_n, j).
\end{equation*}
Furthermore, we also derive identities for $U_{\omega}^{(a)}(\boldsymbol{k})$ with weights being the logarithms, the Gamma function, Bernoulli polynomials and binomials coefficients. 
Our results, which are generalizations of many known results, will be presented in Sections~\ref{sec2}-\ref{sec4}.  The proofs will be given in Section~\ref{sec5}.
\section{The function $U_{\omega}^{(a)}(\boldsymbol{k})$}
\label{sec2}
We recall that an arithmetic function $f$ is called completely multiplicative or additive if $f(mn) =
f(m)f(n)$ or $f(mn)=f(m)+f(n)$, respectively, for all positive integers $m$ and $n$.

In this section, we derive certain identities for $U_{\omega}^{(a)}(\boldsymbol{k})$ when the weight $\omega$ is completely multiplicative or completely additive, respectively. 
\subsection{Completely multiplicative $\omega$}
We are going to state the first main result of the present paper:
\begin{theorem}
\label{Thm1}
Let $(f_i)_{i=1}^n$, $(g_i)_{i=1}^n$ and $(h_i)_{i=1}^n$ be finite sequences of any arithmetic functions and let $\omega$ be a completely multiplicative function. Then we have 
\begin{equation}
\label{eq1}
    U_{\omega}^{(a)}(\boldsymbol{k})
    =\sum_{\boldsymbol{d}|\boldsymbol{k}} \omega^a([\boldsymbol{d}])\boldsymbol{f}(\boldsymbol{d})\boldsymbol{g}\left(\frac{\boldsymbol{k}}{\boldsymbol{d}}\right)
     \sum_{\ell=1}^{[\boldsymbol{k}]^a/[\boldsymbol{d}]^a}\omega(\ell) \boldsymbol{h}\left(\frac{[\boldsymbol{d}]^a \ell}{\boldsymbol{d}^a} \right),
\end{equation}
where, as usual, $\omega^a(x)=(\omega(x))^a$.

In addition, if $h_1, \ldots, h_n$ are completely multiplicative functions, then 
\begin{equation}
\label{eq2}
    U_{\omega}^{(a)}(\boldsymbol{k})
    =\sum_{\boldsymbol{d}|\boldsymbol{k}}\omega^a([\boldsymbol{d}]) \boldsymbol{f}(\boldsymbol{d})\boldsymbol{g}\left(\frac{\boldsymbol{k}}{\boldsymbol{d}}\right) \boldsymbol{h}^a\left(\frac{[\boldsymbol{d}]}{\boldsymbol{d}}\right)
    \sum_{\ell=1}^{[\boldsymbol{k}]^a/[\boldsymbol{d}]^a}\omega(\ell) \boldsymbol{h}(\ell). 
\end{equation}
\end{theorem}

We deduce several corollaries for $\widetilde{U}_{\omega}^{(a)}$ and $\omega=\id_r$, namely: 
\begin{corol}
\label{pro1}
With the above notation we have 
\begin{equation}
\label{eq5}
\widetilde{U}_{\id_r}^{(a)}(\boldsymbol{k})=\frac{[\boldsymbol{k}]^{a r}}{2} \prod_{j=1}^{n}f_j*g_j(k_j)
+\frac{[\boldsymbol{k}]^{a (r+1)}}{r+1}\sum_{m=0}^{\lfloor r/2\rfloor}\binom{r+1}{2m} \frac{B_{2m}}{[\boldsymbol{k}]^{2 a m}}\sum_{\boldsymbol{d}|\boldsymbol{k}}[\boldsymbol{d}]^{a(2m-1)}\boldsymbol{f}(\boldsymbol{d})\boldsymbol{g}\left(\frac{\boldsymbol{k}}{\boldsymbol{d}}\right).
\end{equation}
\end{corol}
From  Eq.~\eqref{eq5}, we get some identities for the weighted average of the products of $\gcd$-sum functions. Taking $f_i*\mu$ in place of $f_i$ and $g_1=\cdots=g_n={\bf 1}$ in Eq.~\eqref{eq5}, we deduce:
\begin{corol}
With the above notation we have 
\begin{multline}
\label{eq8}
\frac{1}{[\boldsymbol{k}]^{ar}}\sum_{j=1}^{[\boldsymbol{k}]^a}j^r\prod_{i=1}^n\Bigl( \sum_{\substack{d_i|k_i\\ d_i^a|j}}f_i*\mu(d)\Bigr)
= \\\frac{\boldsymbol{f}(\boldsymbol{k})}{2}
+\frac{1}{r+1}\sum_{m=0}^{\lfloor r/2\rfloor}\binom{r+1}{2m} \frac{B_{2m}}{[\boldsymbol{k}]^{a(2m-1)}}\sum_{\boldsymbol{d}|\boldsymbol{k}}[\boldsymbol{d}]^{a(2m-1)}\prod_{i=1}^n f_i*\mu (d_i).
\end{multline}
\end{corol}
Taking $a=1$, the left-hand side of Eq.~\eqref{eq8} is then the sum 
$$\frac{1}{[\boldsymbol{k}]^r}\sum_{j=1}^{[\boldsymbol{k}]}j^r \prod_{i=1}^nf_i(\gcd(k_i,j)).$$

Similarly, taking $f_1=\cdots=f_n={\bf 1}$ and $g_i*\mu$ in place of $g_i$ in Eq.~\eqref{eq5}, we deduce:
\begin{corol}
With the above notation we have 
\begin{multline*}
\frac{1}{[\boldsymbol{k}]^{ar}}\sum_{j=1}^{[\boldsymbol{k}]^a}j^r\prod_{i=1}^n\Bigl(\sum_{\substack{d_i|k_i \\ d_i^a|j}} g_i*\mu \left(\frac{k_i}{d_i}\right) \Bigr)
= \frac{\boldsymbol{g}(\boldsymbol{k})}{2}\\
+\frac{1}{r+1}\sum_{m=0}^{\lfloor r/2\rfloor}\binom{r+1}{2m} \frac{B_{2m}}{[\boldsymbol{k}]^{a(2m-1)}}\sum_{\boldsymbol{d}|\boldsymbol{k}}[\boldsymbol{d}]^{a(2m-1)}\prod_{i=1}^ng_i*\mu\left(\frac{k_i}{d_i}\right).
\end{multline*}
\end{corol}
\subsection{Completely additive $\omega$}
In the following theorem we state our second main result:
\begin{theorem}
\label{Thm2}
Let $(f_i)_{i=1}^n$, $(g_i)_{i=1}^n$ and $(h_i)_{i=1}^n$ be finite sequences of any arithmetic functions and let $\omega$ be a completely additive function. Then we have 
\begin{multline}
\label{eq3}
    U_{\omega}^{(a)}(\boldsymbol{k})=\sum_{\boldsymbol{d}|\boldsymbol{k}}\omega\left([\boldsymbol{d}]^a \right)\boldsymbol{f}(\boldsymbol{d})\boldsymbol{g}\left(\frac{\boldsymbol{k}}{\boldsymbol{d}}\right)  
   \sum_{\ell=1}^{[\boldsymbol{k}]^a/[\boldsymbol{d}]^a} \boldsymbol{h}\left(\frac{[\boldsymbol{d}]^a\ell}{\boldsymbol{d}^a} \right)\\
    +
    \sum_{\boldsymbol{d}|\boldsymbol{k}}\boldsymbol{f}(\boldsymbol{d})\boldsymbol{g}\left(\frac{\boldsymbol{k}}{\boldsymbol{d}}\right)  
   \sum_{\ell=1}^{[\boldsymbol{k}]^a/[\boldsymbol{d}]^a} \omega(\ell)\textbf{ h}\left(\frac{[\boldsymbol{d}]^a}{\boldsymbol{d}^a} \ell \right).
\end{multline}
In addition, if $h_1, \ldots, h_n$ are completely multiplicative functions, then 
\begin{multline}
\label{eq4}
    U_{\omega}^{(a)}(\boldsymbol{k})=\sum_{\boldsymbol{d}|\boldsymbol{k}}\omega\left([\boldsymbol{d}]^a\right)\boldsymbol{f}(\boldsymbol{d})\boldsymbol{g}\left(\frac{\boldsymbol{k}}{\boldsymbol{d}}\right) \boldsymbol{h}^a\left(\frac{[\boldsymbol{d}]}{\boldsymbol{d}}\right)
    \sum_{\ell=1}^{[\boldsymbol{k}]^a/[\boldsymbol{d}]^a}\boldsymbol{h}(\ell)\\
    +
    \sum_{\boldsymbol{d}|\boldsymbol{k}}\boldsymbol{f}(\boldsymbol{d})\boldsymbol{g}\left(\frac{\boldsymbol{k}}{\boldsymbol{d}}\right) \boldsymbol{h}^a\left(\frac{[\boldsymbol{d}]}{\boldsymbol{d}}\right)\sum_{\ell=1}^{[\boldsymbol{k}]^a/[\boldsymbol{d}]^a}\omega(\ell)\boldsymbol{h}(\ell). 
\end{multline}
\end{theorem}
Again we deduce several corollaries for $\widetilde{U}_{\omega}^{(a)}$ by taking the weight $\omega=\log$. We have
\begin{corol}
\label{pro2}
With the above notation we have 
\begin{eqnarray}
\label{eq10}
\widetilde{U}_{\log}^{(a)}(\boldsymbol{k})
&=& [\boldsymbol{k}]^a(\log [\boldsymbol{k}]^a-1) \sum_{\boldsymbol{d}|\boldsymbol{k}}\frac{\boldsymbol{f}(\boldsymbol{d})\boldsymbol{g}\left(\frac{\boldsymbol{k}}{\boldsymbol{d}}\right)}{[\boldsymbol{d}]^a} \nonumber
-\frac{a}{2}\sum_{\boldsymbol{d}|\boldsymbol{k}}\log [\boldsymbol{d}] \ \boldsymbol{f}(\boldsymbol{d})\boldsymbol{g}\left(\frac{\boldsymbol{k}}{\boldsymbol{d}}\right) \nonumber\\
&&+\log \sqrt{2\pi [\boldsymbol{k}]^a}\sum_{\boldsymbol{d}|\boldsymbol{k}}\boldsymbol{f}(\boldsymbol{d})\boldsymbol{g}\left(\frac{\boldsymbol{k}}{\boldsymbol{d}}\right)
+\frac{\theta}{12[\boldsymbol{k}]^a}\sum_{\boldsymbol{d}|\boldsymbol{k}} [\boldsymbol{d}]^a\boldsymbol{f}(\boldsymbol{d})\boldsymbol{g}\left(\frac{\boldsymbol{k}}{\boldsymbol{d}}\right),
\end{eqnarray}
for some $\theta \in (0,1)$.
\end{corol}

Taking $g_1=\cdots=g_n={\bf 1}$ and $f_1=\cdots=f_n=\phi$ into Eq.~\eqref{eq10}, we deduce an identity for the weighted average of generalized $\gcd$-sum functions, namely: 
\begin{corol}
With the above notation we have 
\begin{eqnarray*}
\sum_{j=1}^{[\boldsymbol{k}]^a}\log j \prod_{i=1}^n \Bigl(\sum_{d_i|k_i \atop d_i^a|j} \phi(d_i)\Bigr)
&=&  [\boldsymbol{k}]^a(\log [\boldsymbol{k}]^a-1) \sum_{\boldsymbol{d}|\boldsymbol{k}}\frac{\prod_{i=1}^{n}\phi(d_i)}{[\boldsymbol{d}]^a}  
 -\frac{a}{2}\sum_{\boldsymbol{d}|\boldsymbol{k}}\log [\boldsymbol{d}]\prod_{i=1}^{n}\phi(d_i)\\&&+\log \sqrt{2\pi [\boldsymbol{k}]^a}
\prod_{i=1}^{n}\phi(k_i) +\frac{\theta}{12[\boldsymbol{k}]^a}\sum_{\boldsymbol{d}|\boldsymbol{k}}[\boldsymbol{d}]^a\prod_{i=1}^{n}\phi(d_i),  
\end{eqnarray*}
for some $\theta \in (0,1)$.
\end{corol}
When $a=1$, the left-hand side of the above becomes 
$$\sum_{j=1}^{[\boldsymbol{k}]}\log j \prod_{i=1}^n\gcd(k_i, j).$$

\section{Another representation of $\widetilde{U}_{\omega}^{(a)}$}
\label{sec3}
For any two integers $k$ and $j$, the generalized $\gcd$ function $(j, k^a)_a$ is defined as the largest $d\in \mathbb{N}$ such that $d|k$ and $d^a|j$. When $a=1$ the generalized $\gcd$ function becomes the usual $\gcd$ function, see~\cite[Definition 3]{N}. In this section, we use the Dirichlet convolution to give another representation of the function $\widetilde{U}_{\omega}^{(a)}$. 
\begin{theorem}
\label{Thm4}
Let $n\geq 1$ be an integer and $\boldsymbol{k}=(k_1, \ldots, k_n)$ be a vector  with positive integer entries and let $\omega$ be a completely multiplicative function. Then we have 
\begin{equation}\label{eq16}
\widetilde{U}_{\omega}^{(a)}(\boldsymbol{k}) =
\left(\omega \prod_{i=1}^n s^{(a)}_{f_i, g_i, {\bf 1}}\left(k_i,\cdot \right)\right)\ast \Psi^{}([\boldsymbol{k}]^a),
\end{equation}
where 
\begin{equation*}
    \Psi^{}(m)=\sum_{\ell=1 \atop (\ell,m)=1}^{m}\omega (\ell).
\end{equation*}
If $\omega$ is a completely additive function, then we have 
\begin{equation}
\label{eq17}
\widetilde{U}_{\omega}^{(a)}(\boldsymbol{k})=\Bigl(\omega \prod_{i=1}^n s^{(a)}_{f_i, g_i, {\bf 1}}(k_i,\cdot)\Bigr)*\phi([\boldsymbol{k}]^a)+
\Bigl(\prod_{i=1}^n s^{(a)}_{f_i, g_i, {\bf 1}}(k_i,\cdot)\Bigr)*\Psi([\boldsymbol{k}]^a).
\end{equation}
\end{theorem}
For the following we need the Jordan totient function of order $m\in \mathbb{N}$ which is defined as $$\phi_m(n)=n^m\prod_{p \mid n}\left(1-\frac{1}{p^m}\right) = \sum_{d|n}d^m \mu\left(\frac{n}{d}\right).$$

Using Theorem~\ref{Thm4} with $\omega=\id_r$ and $\omega=\log$, respectively, we get: 

\begin{corol}
\label{pro3}
We have 
\begin{multline}
\label{eq18}
\frac{1}{[\boldsymbol{k}]^{ra}}\sum_{j=1}^{[\boldsymbol{k}]^a}j^r \prod_{i=1}^n s_{f_i,g_i,{\bf 1}}^{(a)}(k_i,j)=\frac{1}{2} \prod_{i=1}^n s_{f_i,g_i,{\bf 1}}^{(a)}(k_i,[\boldsymbol{k}]^a)\\
+\frac{1}{r+1}\sum_{m=0}^{\lfloor r/2\rfloor}\binom{r+1}{2m}B_{2m}\sum_{d|[\boldsymbol{k}]^a} d\, \phi_{1-2m}(d)\ \prod_{i=1}^n s_{f_i,g_i,{\bf 1}}^{(a)}\left(k_i,\frac{[\boldsymbol{k}]^a}{d}\right),
\end{multline}
and
\begin{multline}
\label{eq19}
\sum_{j=1}^{[\boldsymbol{k}]^a} \log j\ \prod_{i=1}^n s_{f_i,g_i, {\bf 1}}^{(a)}(k_i,j)=\left(\log \cdot \prod_{i=1}^n s^{(a)}_{f_i, g_i, {\bf 1}}(k_i,\cdot)\right)*\phi([\boldsymbol{k}]^a)\\
 +\sum_{d|[\boldsymbol{k}]^a} \prod_{i=1}^{n}s_{f_i,g_i,{\bf 1}}^{(a)}\left(k_i,\frac{[\boldsymbol{k}]^a}{d}\right)   \left(\sum_{z|d}\mu\left(\frac{d}{z}\right)\log (z!)-\phi(d)\sum_{p|d}\frac{\log p}{p-1} \right),
\end{multline}
where the last sum above is extended over all prime numbers $p$ such that $p|d$. 
\end{corol}
Take $f_1=\cdots=f_n=\phi$ and $g_1=\cdots=g_n={\bf 1}$ in Eq.~\eqref{eq18} and use the identity $\sum_{d|N}\phi(d)=N$ to get: 
\begin{corol}
We have 
\begin{equation*}
\frac{1}{[\boldsymbol{k}]^{ra}}\sum_{j=1}^{[\boldsymbol{k}]^a}j^r \prod_{i=1}^n \Bigl(\sum_{d|k_i \atop d^a|j}\phi(d)\Bigr)
 = \frac{1}{2} \prod_{i=1}^n k_i +\frac{1}{r+1}\sum_{m=0}^{\lfloor r/2\rfloor}\binom{r+1}{2m}B_{2m}\sum_{d|[\boldsymbol{k}]^a} d\, \phi_{1-2m}(d) \prod_{i=1}^n \gcd\left(k_i, \frac{[\boldsymbol{k}]^a}{d}\right). 
\end{equation*}
\end{corol}
\section{Other weighted averages}
\label{sec4}
For $x>0$, the Gamma function $\Gamma$ is defined by 
$$
\Gamma(x)=\int_0^{\infty}\e^{-u}u^{x-1}  \rd u.
$$
We recall that the Bernoulli polynomials are defined by the generating function
\begin{equation*}
\frac{t \e^{x t}}{\e^{t}-1}=\sum_{n=0}^{\infty}B_n(x)\frac{t^n}{n!},
\end{equation*}
for $ |t|<2\pi$. For $x=0$, the numbers $B_n=B_n(0)$  are the Bernoulli numbers. 
In this section, we consider weighted averages of the product $s_{f_1,g_1,h_1}^{(a)}(k,j)\cdots s_{f_n, g_n,h_n}^{(a)}(k,j)$ with weights being the Gamma function, Binomial coefficients, and Bernoulli polynomials. We prove that: 
\begin{theorem}
\label{Thm3}
Let $n\geq 1$ be an integer and $\boldsymbol{k}=(k_1, \ldots, k_n)$ be a vector  with positive integer entries and let $a\ge 2$ be a fixed integer. Then we have the following formulas:
\begin{multline}
\label{eq13}
\sum_{j=1}^{[\boldsymbol{k}]^a}\log\Gamma \left(\frac{j}{[\boldsymbol{k}]^a}\right)\ \prod_{i=1}^n s_{f_i,g_i, {\bf 1}}^{(a)}(k_i,j)
 =  \frac{[\boldsymbol{k}]^a}{2}\log(2\pi) \sum_{\boldsymbol{d}|\boldsymbol{k}}\frac{\boldsymbol{f}(\boldsymbol{d})\boldsymbol{g}\left(\frac{\boldsymbol{k}}{\boldsymbol{d}}\right)}{[\boldsymbol{d}]^a}
\\-\frac{1}{2}\log(2\pi [\boldsymbol{k}]^a)\prod_{i=1}^{n}f_i*g_i(k_i) 
+\frac{a}{2} \sum_{\boldsymbol{d}|\boldsymbol{k}}\log[\boldsymbol{d}]\ \boldsymbol{f}(\boldsymbol{d})\boldsymbol{g}\left(\frac{\boldsymbol{k}}{\boldsymbol{d}}\right),
\end{multline}
\begin{equation}
\label{eq14}
\sum_{j=0}^{[\boldsymbol{k}]^a}\binom{[\boldsymbol{k}]^a}{j}\ \prod_{i=1}^n s_{f_i,g_i, {\bf 1}}^{(a)}(k_i,j)
=2^{[\boldsymbol{k}]^a}\sum_{\boldsymbol{d}|\boldsymbol{k}}\frac{\boldsymbol{f}(\boldsymbol{d})\boldsymbol{g}\left(\frac{\boldsymbol{k}}{\boldsymbol{d}}\right)}{[\boldsymbol{d}]^a}\sum_{\ell=1}^{[\boldsymbol{d}]^a}(-1)^{ [\boldsymbol{k}]^a\ell/[\boldsymbol{d}]^a}\cos^{[\boldsymbol{k}]^a}\left( \frac{\pi \ell }{[\boldsymbol{d}]^a}\right),
\end{equation}
and 
\begin{equation}
\label{eq15}
\sum_{j=0}^{[\boldsymbol{k}]^a-1} B_m\left(\frac{j}{[\boldsymbol{k}]^a}\right)\ \prod_{i=1}^n s_{f_i,g_i, {\bf 1}}^{(a)}(k_i,j)= \frac{B_m}{[\boldsymbol{k}]^{a(m-1)}}\sum_{\boldsymbol{d}|\boldsymbol{k}}\frac{\boldsymbol{f}(\boldsymbol{d})\boldsymbol{g}\left(\frac{\boldsymbol{k}}{\boldsymbol{d}}\right)}{[\boldsymbol{d}]^{a(1-m)}}.
\end{equation}
\end{theorem}
Now we provide an application of Theorem~\ref{Thm3}. In particular, we derive formulas for weighted averages of the product of the $\gcd$-sum function with weights being the Gamma function, Binomial coefficients, and Bernoulli polynomials by taking $f_i=\phi$ and $g_i={\bf 1}$ for all $i= 1,\ldots, n$ into Eqs.~\eqref{eq13}, \eqref{eq14} and \eqref{eq15}. We obtain the following results: 
\begin{corol}
Under the hypothesis of Theorem~\ref{Thm3}, we have 
\begin{eqnarray*}
\sum_{j=1}^{[\boldsymbol{k}]^a}\log\Gamma\Bigl(\frac{j}{[\boldsymbol{k}]^a}\Bigr)  \prod_{i=1}^{n}\Bigg( \sum_{d_i|k_i\atop d_i^{a}|j} \phi(d_i)\Bigg)
&=& \frac{[\boldsymbol{k}]^a}{2}\log(2\pi) \sum_{\boldsymbol{d}|\boldsymbol{k}}\frac{\prod_{i=1}^{n}\phi(d_i)}{[\boldsymbol{d}]^a} \\&&-\frac{k_1\cdots k_n}{2}\log(2\pi [\boldsymbol{k}]^a) +\frac{a}{2}\sum_{\boldsymbol{d}|\boldsymbol{k}}\log [\boldsymbol{d}]\ \prod_{i=1}^{n}\phi(d_i),
\end{eqnarray*}
\begin{equation*}
\sum_{j=0}^{[\boldsymbol{k}]^a}\binom{[\boldsymbol{k}]^a}{j}\prod_{i=1}^{n}\Bigg( \sum_{d_i|k_i\atop d_i^{a}|j} \phi(d_i)\Bigg)
=2^{[\boldsymbol{k}]^a}\sum_{\boldsymbol{d}|\boldsymbol{k}}\frac{\phi(d_1)\cdots \phi(d_n)}{[\boldsymbol{d}]^a}\sum_{\ell=1}^{[\boldsymbol{d}]^a}(-1)^{\frac{\ell [\boldsymbol{k}]^a}{[\boldsymbol{d}]^a}}\cos^{[\boldsymbol{k}]^a}\left(\frac{\pi \ell}{ [\boldsymbol{d}]^a}\right),
\end{equation*}
and 
$$
\sum_{j=0}^{[\boldsymbol{k}]^a-1}B_m\left(\frac{j}{[\boldsymbol{k}]^a}\right)\prod_{i=1}^{n}\Bigg( \sum_{d_i|k_i\atop d_i^{a}|j} \phi(d_i)\Bigg)=
\frac{B_m}{[\boldsymbol{k}]^{a(m-1)}}\sum_{\boldsymbol{d}|\boldsymbol{k}}\frac{\phi(d_1)\cdots \phi(d_n)}{[\boldsymbol{d}]^{a(1-m)}}
$$
In the case $m=1$, we have
\begin{equation*}
\sum_{j=0}^{[\boldsymbol{k}]^a-1}B_1\left(\frac{j}{[\boldsymbol{k}]^a}\right)\prod_{i=1}^{n}\Bigg( \sum_{d_i|k_i\atop d_i^{a}|j} \phi(d_i)\Bigg)=
-\frac{k_1\cdots k_n}{2}.
\end{equation*}
\end{corol}

\section{Proofs}
\label{sec5}
\begin{proof}[Proof of Theorem \ref{Thm1}]
We recall that 
\begin{equation*}
s^{(a)}_{f,g,h}(k,j)=\sum_{\substack{d|k\\d^a|j}}f(d)g\left(\frac{k}{d}\right)h\left(\frac{j}{d^a}\right).
\end{equation*} 
Then we get with Eq.~\eqref{defU} 
\begin{eqnarray}
\label{main}
U_{\omega}^{(a)}(\boldsymbol{k})&=& 
\sum_{j=1}^{[\boldsymbol{k}]^a}\omega(j)\sum_{\substack{\boldsymbol{d}|\boldsymbol{k}\\\boldsymbol{d}^a|j}}\boldsymbol{f}(\boldsymbol{d})\boldsymbol{g}\left(\frac{\boldsymbol{k}}{\boldsymbol{d}}\right)\boldsymbol{h}\left(\frac{j}{\boldsymbol{d}^a}\right)\nonumber
\\&=&
\sum_{\boldsymbol{d}|\boldsymbol{k}}\boldsymbol{f}(\boldsymbol{d})\boldsymbol{g}\left(\frac{\boldsymbol{k}}{\boldsymbol{d}}\right)\sum_{\substack{j=1\\\boldsymbol{d}^{a}|j}}^{[\boldsymbol{k}]^a}\omega(j)\boldsymbol{h}\left(\frac{j}{\boldsymbol{d}^a}\right).
\end{eqnarray}
Using the fact that $\omega$ is a completely multiplicative function, we obtain 
\begin{equation*}
    U_{\omega}^{(a)}(\boldsymbol{k})
    =\sum_{\boldsymbol{d}|\boldsymbol{k}}\omega^a([\boldsymbol{d}])\boldsymbol{f}(\boldsymbol{d})\boldsymbol{g}\left(\frac{\boldsymbol{k}}{\boldsymbol{d}}\right)
     \sum_{\ell=1}^{[\boldsymbol{k}]^a/[\boldsymbol{d}]^a}\omega(\ell) \boldsymbol{h}\Bigg(\frac{[\boldsymbol{d}]^a \ell }{\boldsymbol{d}^a}\Bigg). 
\end{equation*}
This proves Eq.~\eqref{eq1}. 

In case that all $h_i$, $i=1, \ldots, n$, are completely multiplicative, the latter function can be rewritten as 
\begin{equation*}
    U_{\omega}^{(a)}(\boldsymbol{k})
    =\sum_{\boldsymbol{d}|\boldsymbol{k}}\omega^a([\boldsymbol{d}])\boldsymbol{f}(\boldsymbol{d})\boldsymbol{g}\left(\frac{\boldsymbol{k}}{\boldsymbol{d}}\right) 
     \boldsymbol{h}^a\left(\frac{[\boldsymbol{d}]}{\boldsymbol{d}}\right)
     \sum_{\ell=1}^{[\boldsymbol{k}]^a/[\boldsymbol{d}]^a}\omega(\ell) \boldsymbol{h}(\ell).
\end{equation*}
This completes the proof of Theorem~\ref{Thm1}. 
\end{proof}
Theorem~\ref{Thm2} can be proven in the same way.
\begin{proof}[Proof of Corollary~\ref{pro1}]
We recall that 
\begin{equation*}
s_{f, g,{\bf 1}}^{(a)}(k,j)=\sum_{\substack{d|k\\d^a|j}}f(d)g\left(\frac{k}{d}\right).
\end{equation*}
Substituting $\omega=\id_r$ and $h_1=\ldots=h_n={\bf 1}$ into Eq.~\eqref{eq2}, we get
\begin{eqnarray*}
\widetilde{U}_{\id_r}^{(a)}(\boldsymbol{k})&=&
\sum_{j=1}^{[\boldsymbol{k}]^a}j^r\sum_{\substack{\boldsymbol{d}|\boldsymbol{k}\\\boldsymbol{d}^a|j}}\boldsymbol{f}(\boldsymbol{d})\boldsymbol{g}\left(\frac{\boldsymbol{k}}{\boldsymbol{d}}\right)
\\&=&
\sum_{\boldsymbol{d}|\boldsymbol{k}}\boldsymbol{f}(\boldsymbol{d})\boldsymbol{g}\left(\frac{\boldsymbol{k}}{\boldsymbol{d}}\right)\sum_{\substack{j=1\\\boldsymbol{d}^a|j}}^{[\boldsymbol{k}]^a}j^r
\\&=&
\sum_{\boldsymbol{d}|\boldsymbol{k}}\boldsymbol{f}(\boldsymbol{d})\boldsymbol{g}\left(\frac{\boldsymbol{k}}{\boldsymbol{d}}\right)[\boldsymbol{d}]^{ar}\sum_{\ell=1}^{[\boldsymbol{k}]^a/[\boldsymbol{d}]^a} \ell^r.
\end{eqnarray*}
 Using the fact that 
\begin{equation*}
\sum_{m=1}^{N}m^r=\frac{N^r}{2}+\frac{1}{r+1}\sum_{m=0}^{\lfloor r/2\rfloor}\binom{r+1}{2m}B_{2m}N^{r+1-2m},
\end{equation*}
for any integer $N\ge 2$, see~\cite[Proposition~9.2.12]{Coh}  or \cite[Section~3.9]{Co}, we deduce that 
\begin{equation*}
\widetilde{U}_{\id_r}^{(a)}(\boldsymbol{k})
  =
  \frac{[\boldsymbol{k}]^{ra}}{2} \prod_{i=1}^n f_i*g_i(k_i)
  +\frac{[\boldsymbol{k}]^{a(r+1)}}{r+1}\sum_{m=0}^{\lfloor r/2\rfloor}\binom{r+1}{2m}\frac{B_{2m}}{[\boldsymbol{k}]^{2ma}}\sum_{\boldsymbol{d}|\boldsymbol{k}}\boldsymbol{f}(\boldsymbol{d})\boldsymbol{g}\left(\frac{\boldsymbol{k}}{\boldsymbol{d}}\right)[\boldsymbol{d}]^{a(2m-1)}.
\end{equation*}
This completes the proof. 
\end{proof}


\begin{proof}[Proof of Corollary~\ref{pro2}]
Substituting $\omega=\log$ and $h_1=\cdots=h_n={\bf 1}$ into Eq.~\eqref{eq3}, we obtain 
\begin{eqnarray*}
\widetilde{U}_{\log}^{(a)}(\boldsymbol{k}) & = & [\boldsymbol{k}]^a\sum_{\boldsymbol{d}|\boldsymbol{k}}\boldsymbol{f}(\boldsymbol{d})\boldsymbol{g}\left(\frac{\boldsymbol{k}}{\boldsymbol{d}}\right)\frac{\log[\boldsymbol{d}]^a}{[\boldsymbol{d}]^a} + \sum_{\boldsymbol{d}|\boldsymbol{k}}\boldsymbol{f}(\boldsymbol{d})\boldsymbol{g}\left(\frac{\boldsymbol{k}}{\boldsymbol{d}}\right)\sum_{\ell=1}^{[\boldsymbol{k}]^a/[\boldsymbol{d}]^a}\log \ell\\
& = &  
    a[\boldsymbol{k}]^a \sum_{\boldsymbol{d}|\boldsymbol{k}}\boldsymbol{f}(\boldsymbol{d})\boldsymbol{g}\left(\frac{\boldsymbol{k}}{\boldsymbol{d}}\right)\frac{\log [\boldsymbol{d}]}{[\boldsymbol{d}]^a}
 +\sum_{\boldsymbol{d}|\boldsymbol{k}}\boldsymbol{f}(\boldsymbol{d})\boldsymbol{g}\left(\frac{\boldsymbol{k}}{\boldsymbol{d}}\right)\log \left(\frac{[\boldsymbol{k}]^a}{[\boldsymbol{d}]^a}\right)!.
\end{eqnarray*}
Using the Stirling formula, see~\cite[p. 91]{Ni},
\begin{equation*}
   \log (\ell !) =\ell \log \ell -\ell+\frac{1}{2}\log \ell+ \log \sqrt{2\pi}+\frac{\theta}{12\ell},
\end{equation*}
the formula~\eqref{eq10} is proved.
\end{proof}

\begin{proof}[Proof of Theorem~\ref{Thm4}]
Since $[\boldsymbol{k}]=\lcm(k_1, \ldots, k_n)$, we can write  
\begin{equation*}
s^{(a)}_{f, g, {\bf 1}}(k_i,j) = \sum_{d|k_i \atop d^a|j, d|[\boldsymbol{k}]} f(d)g\left( \frac{k_i}{d}\right) =  \sum_{d|k_i\atop d^a|(j, [\boldsymbol{k}]^a)} f(d)g\left( \frac{k_i}{d}\right) = s^{(a)}_{f,g, {\bf 1}}\left(k_i,(j, [\boldsymbol{k}]^a)\right).
\end{equation*}
Hence
$$\widetilde{U}_{\omega}^{(a)}(\boldsymbol{k}) = \sum_{j=1}^{[k]^a}\omega (j) \prod_{i=1}^n s^{(a)}_{f_i, g_i, {\bf 1}}\left(k_i,(j, [\boldsymbol{k}]^a)\right).$$
Now we split the range of summation  of the first sum over $j$ according to the value of $(j,[\boldsymbol{k}]^a)$. This gives 
\begin{eqnarray*}
\widetilde{U}_{\omega}^{(a)}(\boldsymbol{k}) & = &  \sum_{d|[\boldsymbol{k}]^a} \prod_{i=1}^n s^{(a)}_{f_i, g_i, {\bf 1}}\left(k_i,d\right) \sum_{j=1 \atop (j,[\boldsymbol{k}]^a)=d}^{[\boldsymbol{k}]^a}\omega (j) \\
& = & \sum_{d|[\boldsymbol{k}]^a} \prod_{i=1}^n s^{(a)}_{f_i, g_i, {\bf 1}}\left(k_i,d\right) \sum_{\ell=1 \atop (\ell,[\boldsymbol{k}]^a/d)=1}^{[\boldsymbol{k}]^a/d}\omega (\ell d).
\end{eqnarray*}
Now we use the fact that $\omega$ is completely multiplicative. This way we obtain
\begin{eqnarray*}
\widetilde{U}_{\omega}^{(a)}(\boldsymbol{k}) & = & \sum_{d|[\boldsymbol{k}]^a} \omega(d) \prod_{i=1}^n s^{(a)}_{f_i, g_i, {\bf 1}}\left(k_i,d\right) \sum_{\ell=1 \atop (\ell,[\boldsymbol{k}]^a/d)=1}^{[\boldsymbol{k}]^a/d}\omega (\ell)\\
& = & \sum_{d|[\boldsymbol{k}]^a} \omega(d) \prod_{i=1}^n s^{(a)}_{f_i, g_i, {\bf 1}}\left(k_i,d\right)\Psi^{}\left(\frac{[\boldsymbol{k}]^a}{d}\right)\\
& = & \Bigg(\omega \prod_{i=1}^n s^{(a)}_{f_i, g_i, {\bf 1}}\left(k_i,\cdot \right)\Bigg)\ast \Psi^{}([\boldsymbol{k}]^a).
\end{eqnarray*}
This proves Eq.~\eqref{eq16}. Eq.~\eqref{eq17} is shown in the same way. 
\end{proof}
\begin{remark}
With the same method we can also prove the following representations: if $\omega$ is completely multiplicative, then
\begin{equation*}
\widetilde{U}_{\omega}^{(a)}(\boldsymbol{k})= \left(\omega^a \prod_{i=1}^n s^{(1)}_{f_i, g_i, {\bf 1}}(k_i,\cdot)\right)*\Psi^{(a)}([\boldsymbol{k}])
\end{equation*}
with $\omega^a(d)=(\omega(d))^a$ and for $N \in \mathbb{N}$, 
\begin{equation*}
    \Psi^{(a)}(N)=\sum_{\ell=1 \atop (\ell, N^a)_a=1}^{N^a}\omega(\ell).
\end{equation*}
If $\omega$ is a completely additive function, then we have 
\begin{equation*}
\widetilde{U}_{\omega}^{(a)}(\boldsymbol{k})=\left(\omega^{(a)} \prod_{i=1}^n s^{(1)}_{f_i, g_i, {\bf 1}}(k_i,\cdot)\right)*\Phi^{(a)}([\boldsymbol{k}])+
\left(\prod_{i=1}^n s^{(1)}_{f_i, g_i, {\bf 1}}(k_i,\cdot)\right)*\Psi^{(a)}([\boldsymbol{k}]),
\end{equation*}
with $\omega^{(a)}(d)=\omega(d^a)$ and for $N \in \mathbb{N}$, $$\Phi^{(a)}(N)=\sum_{\ell=1 \atop (\ell, N^a)_a=1}^{N^a}1=\#\{\ell \in \{1,\ldots,N^a\}\ : \ (\ell, N^a)_a=1\} .$$
\end{remark}
\begin{proof}[Proof of Corollary~\ref{pro3}]
Substituting $\omega=\id_r$ into Eq.~\eqref{eq16},  we get 
\begin{equation*}
    \sum_{j=1}^{[\boldsymbol{k}]^a}j^r \prod_{i=1}^n s_{f_i, g_i, {\bf 1}}^{(a)}(k_i,j)=\sum_{d|[\boldsymbol{k}]^a}d^ r \prod_{i=1}^n s_{f_i, g_i, {\bf 1}}^{(a)}(k_i,d) \sum_{\ell=1 \atop (\ell, [\boldsymbol{k}]^a/d)=1}^{[\boldsymbol{k}]^a/d}\ell^r.
\end{equation*}
Applying the following formula, see~\cite[Corollary 4]{Sin}, to the last sum above
\begin{equation*}
 \sum_{m=1 \atop (m, N)=1}^N m^r=\frac{N^{r+1}}{r+1}\sum_{m=0}^{\lfloor r/2\rfloor}\binom{r+1}{2m}B_{2m}\phi_{1-2m}(N),   
\end{equation*}
for any positive integer $N>1$, we find that 
\begin{eqnarray*}
\lefteqn{\sum_{j=1}^{[\boldsymbol{k}]^a}j^r \prod_{i=1}^{n}s_{f_i,g_i,{\bf 1}}^{(a)}(k_i,j)}\\
&=&
[\boldsymbol{k}]^{ra}\prod_{i=1}^{n}s_{f_i,g_i,{\bf 1}}^{(a)}(k_i,[\boldsymbol{k}]^a)
+\frac{[\boldsymbol{k}]^{ra}}{r+1}\sum_{d|[\boldsymbol{k}]^a\atop d<[\boldsymbol{k}]^a}\prod_{i=1}^{n}s_{f_i,g_i,{\bf 1}}^{(a)}(k_i,d)\sum_{m=0}^{\lfloor r/2\rfloor}\binom{r+1}{2m}B_{2m}\phi_{1-2m}\left(\frac{[\boldsymbol{k}]^a}{d}\right)\frac{[\boldsymbol{k}]^a}{d}
\\&=&
[\boldsymbol{k}]^{ra}\prod_{i=1}^{n}s_{f_i,g_i,{\bf 1}}^{(a)}(k_i,[\boldsymbol{k}]^a)
+\frac{[\boldsymbol{k}]^{ra}}{r+1}\sum_{m=0}^{\lfloor r/2\rfloor}\binom{r+1}{2m}B_{2m}\sum_{\ell|[\boldsymbol{k}]^a\atop \ell>1} \ell \phi_{1-2m}(\ell) \prod_{i=1}^{n}s_{f_i,g_i,{\bf 1}}^{(a)}\left(k_i,\frac{[\boldsymbol{k}]^a}{\ell}\right)
\\&=&
[\boldsymbol{k}]^{ra}\prod_{i=1}^{n}s_{f_i,g_i,{\bf 1}}^{(a)}(k_i,[\boldsymbol{k}]^a)-\frac{[\boldsymbol{k}]^{ra}}{r+1}\prod_{i=1}^{n}s_{f_i,g_i,{\bf 1}}^{(a)}\left(k_i,[\boldsymbol{k}]^a\right)\sum_{m=0}^{\lfloor r/2\rfloor}\binom{r+1}{2m}B_{2m}
\\&&+\frac{[\boldsymbol{k}]^{ra}}{r+1}\sum_{m=0}^{\lfloor r/2\rfloor}\binom{r+1}{2m}B_{2m}\sum_{\ell|[\boldsymbol{k}]^a} \ell \phi_{1-2m}(\ell) \prod_{i=1}^{n}s_{f_i,g_i,{\bf 1}}^{(a)}\left(k_i,\frac{[\boldsymbol{k}]^a}{\ell}\right).
\end{eqnarray*}
Using the well-known formula, 
\begin{equation*}
    \sum_{m=0}^{\lfloor r/2\rfloor}\binom{r+1}{2m}B_{2m}=\frac{r+1}{2},
\end{equation*}
we obtain the desired formula \eqref{eq18}. 

In order to prove Eq.~\eqref{eq19} we substitute $\omega=\log$ into Eq.~\eqref{eq17} and use the formula, see~\cite[Lemma 12]{I.K.M},
\begin{equation*}
\sum_{\ell=1\atop (\ell, N)=1}^N\log \ell=\sum_{d|N}\mu\left(\frac{N}{d}\right)\log (d!)-\phi(N)\sum_{p|N}\frac{\log p}{p-1},
\end{equation*}
for any positive integer $N$, where the last sum above is extended over all prime numbers $p$ such that $p|N$. This gives 
\begin{eqnarray*}
\sum_{j=1}^{[\boldsymbol{k}]^a} \log j \prod_{i=1}^n s_{f_i,g_i,{\bf 1}}^{(a)}(k_i,j)
& = & \Bigl(\log \cdot \prod_{i=1}^n s^{(a)}_{f_i, g_i, {\bf 1}}(k_i,\cdot)\Bigr)*\phi([\boldsymbol{k}]^a)+
\Bigl(\prod_{i=1}^n s^{(a)}_{f_i, g_i, {\bf 1}}(k_i,\cdot)\Bigr)*\Psi([\boldsymbol{k}]^a)\\
& = &\Bigl(\log \cdot \prod_{i=1}^n s^{(a)}_{f_i, g_i, {\bf 1}}(k_i,\cdot)\Bigr)*\phi([\boldsymbol{k}]^a)+
\sum_{d|[\boldsymbol{k}]^a} \prod_{i=1}^{n}s_{f_i,g_i,{\bf 1}}^{(a)}\left(k_i,\frac{[\boldsymbol{k}]^a}{d}\right)\sum_{\ell=1 \atop (\ell,d)=1}^{d} \log \ell\\
& = & \Bigl(\log \cdot \prod_{i=1}^n s^{(a)}_{f_i, g_i, {\bf 1}}(k_i,\cdot)\Bigr)*\phi([\boldsymbol{k}]^a)\\
& & +\sum_{d|[\boldsymbol{k}]^a} \prod_{i=1}^{n}s_{f_i,g_i,{\bf 1}}^{(a)}\left(k_i,\frac{[\boldsymbol{k}]^a}{d}\right)   \Bigg(\sum_{z|d}\mu\left(\frac{d}{z}\right)\log (z!)-\phi(d)\sum_{p|d}\frac{\log p}{p-1} \Bigg)
\end{eqnarray*}
as desired.
\end{proof}
\begin{proof}[Proof of Theorem~\ref{Thm3}]
Substituting $\omega(d)=\log \Gamma(d/[\boldsymbol{k}]^a)$ and $h_1=\ldots=h_n={\bf 1}$ into Eq.~\eqref{main}, we get 
\begin{equation*}
\sum_{j=1}^{[\boldsymbol{k}]^a}\log \Gamma\left(\frac{j}{[\boldsymbol{k}]^a}\right)\prod_{i=1}^ns_{f_i,g_i,{\bf 1}}^{(a)}(k_i,j)
=\sum_{\boldsymbol{d}|\boldsymbol{k}}\boldsymbol{f}(\boldsymbol{d})\boldsymbol{g}\left(\frac{\boldsymbol{k}}{\boldsymbol{d}}\right)\sum_{\ell=1}^{[\boldsymbol{k}]^a/[\boldsymbol{d}]^a}\log \Gamma\left(\frac{\ell}{[\boldsymbol{k}]^a/[\boldsymbol{d}]^a}\right).
\end{equation*}
From the Gauss-Legendre multiplication formula for the Gamma function it is known that 
\begin{equation*}
\prod_{j=1}^{n}\Gamma\left( \frac{j}{n}\right)=\frac{(2\pi)^{\frac{n-1}{2}}}{\sqrt{n}} \ \ \ \mbox{ for all $n \in \mathbb{N}$.}
\end{equation*}
Hence it follows that 
\begin{eqnarray*}
\sum_{j=1}^{[\boldsymbol{k}]^a}\log\Gamma\left(\frac{j}{[\boldsymbol{k}]^a}\right)\prod_{i=1}^ns_{f_i,g_i,{\bf 1}}^{(a)}(k_i,j)
& =& \sum_{\boldsymbol{d}|\boldsymbol{k}}\boldsymbol{f}(\boldsymbol{d})\boldsymbol{g}\left(\frac{\boldsymbol{k}}{\boldsymbol{d}}\right)\left[\frac{[\boldsymbol{k}]^a}{2[\boldsymbol{d}]^a} \log(2\pi)- \frac{1}{2} \log(2\pi [\boldsymbol{k}]^a)+ \frac{a}{2}\log [\boldsymbol{d}]\right]\\
& =&\frac{[\boldsymbol{k}]^a}{2}\log(2\pi) \sum_{\boldsymbol{d}|\boldsymbol{k}}\frac{\boldsymbol{f}(\boldsymbol{d})\boldsymbol{g}(\boldsymbol{k}/\boldsymbol{d})}{[\boldsymbol{d}]^a} -\frac{1}{2} \log(2\pi [\boldsymbol{k}]^a) \prod_{i=1}^n f_i\ast g_i(k_i) \\& & + \frac{a}{2} \sum_{\boldsymbol{d}|\boldsymbol{k}}  \log [\boldsymbol{d}]\ \boldsymbol{f}(\boldsymbol{d})\boldsymbol{g}\left(\frac{\boldsymbol{k}}{\boldsymbol{d}}\right). 
\end{eqnarray*}
This is precisely Eq.~\eqref{eq13}. 

Now set $w(j)=\binom{K^a}{j}$ and $h_1=\ldots=h_n={\bf 1}$ into Eq.~\eqref{main} to get

\begin{equation*}
\sum_{j=0}^{[\boldsymbol{k}]^a}\binom{[\boldsymbol{k}]^a}{j}\prod_{i=1}^ns_{f_i,g_i,{\bf 1}}^{(a)}(k_i,j)=\sum_{\boldsymbol{d}|\boldsymbol{k}}\boldsymbol{f}(\boldsymbol{d})\boldsymbol{g}\left(\frac{\boldsymbol{k}}{\boldsymbol{d}}\right)\sum_{m=0}^{[\boldsymbol{k}]^a/[\boldsymbol{d}]^a}\binom{[\boldsymbol{k}]^a}{m[\boldsymbol{d}]^a}.
\end{equation*}

For any positive integers $r$ and $n$ it is known that
\begin{equation*}
    \sum_{m=0}^{\lfloor n/r\rfloor}\binom{n}{mr}=\frac{2^n}{r}\sum_{\ell=1}^{r}\cos^n  \left(\frac{\pi \ell}{r}\right)\cos \left(\frac{\pi \ell n}{r}\right),
\end{equation*}
see \cite[Eq.~(27)]{To3}. Hence 
\begin{eqnarray*}
\sum_{j=0}^{[\boldsymbol{k}]^a}\binom{[\boldsymbol{k}]^a}{j}\prod_{i=1}^ns_{f_i,g_i,{\bf 1}}^{(a)}(k_i,j)
&= & 2^{[\boldsymbol{k}]^a}\sum_{\boldsymbol{d}|\boldsymbol{k}}\frac{\boldsymbol{f}(\boldsymbol{d})\boldsymbol{g}\left(\boldsymbol{k}/ \boldsymbol{d}\right)}{[\boldsymbol{d}]^a}\sum_{\ell=1}^{[\boldsymbol{d}]^a}\cos ^{[\boldsymbol{k}]^a}\left(\frac{\pi \ell}{[\boldsymbol{d}]^a}\right)\cos \left( \frac{\pi \ell [\boldsymbol{k}]^a}{[\boldsymbol{d}]^a}\right)\\
& = & 
2^{[\boldsymbol{k}]^a}\sum_{\boldsymbol{d}|\boldsymbol{k}}\frac{\boldsymbol{f}(\boldsymbol{d})\boldsymbol{g}\left(\boldsymbol{k}/\boldsymbol{d}\right)}{[\boldsymbol{d}]^a}\sum_{\ell=1}^{[\boldsymbol{d}]^a}(-1)^{\ell [\boldsymbol{k}]^a/[\boldsymbol{d}]^a}\cos ^{[\boldsymbol{k}]^a}\left(\frac{\pi \ell}{[\boldsymbol{d}]^a}\right).
\end{eqnarray*}
This completes the proof of Eq.~\eqref{eq14}. 

Finally, put $w(j)=B_m(j/[\boldsymbol{k}]^a)$ and $h_1=\ldots=h_n={\bf 1}$ in Eq.~\eqref{main}. Then  
\begin{equation*}
\sum_{j=0}^{[\boldsymbol{k}]^a-1}B_m\left(\frac{j}{[\boldsymbol{k}]^a}\right)\prod_{i=1}^ns_{f_i,g_i,{\bf 1}}^{(a)}(k_i,j)=\sum_{\boldsymbol{d}|\boldsymbol{k}}\boldsymbol{f}(\boldsymbol{d})\boldsymbol{g}\left(\frac{\boldsymbol{k}}{\boldsymbol{d}}\right)\sum_{\ell=0}^{[\boldsymbol{k}]^a/[\boldsymbol{d}]^a-1}B_m\left(\frac{[\boldsymbol{d}]^a\ell}{[\boldsymbol{k}]^a}\right).
\end{equation*}
Using the following property of the Bernoulli polynomials, see~\cite[Proposition 9.1.3]{Coh}, 
\begin{equation*}
    \sum_{\ell=0}^{k-1}B_m\left(\frac{\ell}{k}\right)=\frac{B_m}{k^{m-1}},
\end{equation*}
we find that 
\begin{equation*}
\sum_{j=0}^{[\boldsymbol{k}]^a-1}B_m\left(\frac{j}{[\boldsymbol{k}]^a}\right)\prod_{i=1}^ns_{f_i,g_i,{\bf 1}}^{(a)}(k_i,j)=
\frac{B_m}{[\boldsymbol{k}]^{a(m-1)}}
\sum_{\boldsymbol{d}|\boldsymbol{k}}[\boldsymbol{d}]^{a(m-1)}\boldsymbol{f}(\boldsymbol{d})\boldsymbol{g}\left(\frac{\boldsymbol{k}}{\boldsymbol{d}}\right).
\end{equation*}
This completes the proof of Theorem~\ref{Thm3}.
\end{proof}



\medskip\noindent {\footnotesize Isao Kiuchi: 
Department of Mathematical Sciences,
Faculty of Science, Yamaguchi University,
Yoshida 1677-1, Yamaguchi 753-8512, Japan.
e-mail: {\tt kiuchi@yamaguchi-u.ac.jp}}

\medskip\noindent {\footnotesize Friedrich Pillichshammer: 
Institute of Financial Mathematics and Applied Number Theory,
Johannes Kepler University,
Altenbergerstrasse 69, 4040 Linz, Austria.
e-email: {\tt friedrich.pillichshammer@jku.at}}

\medskip\noindent {\footnotesize Sumaia Saad Eddin: 
Institute of Financial Mathematics and Applied Number Theory,
Johannes Kepler University,
Altenbergerstrasse 69, 4040 Linz, Austria.
e-mail: {\tt sumaia:saad\_eddin@jku.at}}


\begin{thebibliography}{9}

\bibitem{AA} D. R.  Anderson and T. M. Apostol,  The evaluation of Ramanujan's sum and generalizations, {\it Duke Math. J.} {\bf 20}  (1953),  211--216.


\bibitem{C1} E. Cohen, An extension of Ramanujan's sum, 
 {\it Duke Math. J.} {\bf 16} (1949), 85--90.

\bibitem{C2} E. Cohen, An extension of Ramanujan's sum. II. Additive properties, {\it Duke Math. J.} {\bf 22} (1955), 543--550.


\bibitem{C3} E. Cohen, An extension of Ramanujan's sum. III. Connections with totient functions), {\it Duke Math. J.} {\bf 23} (1956), 623--630.


\bibitem{Coh} H. Cohen, {\it Number Theory, Volume II: Analytic and Modern Tools, Graduate Texts in Mathematics},  {\bf 240}, Springer, 2007.

\bibitem{Co} L. Comtet, {\it Advanced Combinatorics, The Art of Finite and Infinite Expansions}, D. Reidel Publishing Co, Springer, 1974.


\bibitem{I.K.M} S. Ikeda, I. Kiuchi and K. Matsuoka, Sums of products of generalized Ramanujan sums, {\it J. Integer Sequences} {\bf 19} (2016), Article 16.2.7.

\bibitem{K2017} I. Kiuchi, Sums of averages of generalized Ramanujan sums,  {\it J. Number Theory} {\bf 180} (2017), 310-348.



\bibitem{K.R1} P. K\"uhn and N. Robles, Explicit formulas of a generalized Ramanujan sum,  {\it Int. J. Number Theory} {\bf 12} (2016), 383-408.

\bibitem{L1} V. A. Liskovets, A multivariate arithmetic function of combinatorial and topological significance, {\it Integers} \textbf{10: A12} (2010), 155-177.


\bibitem{Mc1} P. J. McCarthy, {\it Introduction to Arithmetical Functions}, Springer, 1986.


\bibitem{Mc3} P. J. McCarthy, The generation of arithmetical identities, {\it J. Reine Angew. Math.} {\bf 203} (1960), 55-63. 


\bibitem{Mc2} P. J. McCarthy, Some properties of the extended Ramanujan sums, {\it Arch. der Math.} {\bf 11} (1960), 253-258. 


\bibitem{M.N1} A. D. Mednykh and R. Nedela, Enumeration of unrooted maps of a given genus, {\it J. Combin. Theory Ser. B} {\bf 96} (2006), 706-729.

\bibitem{N} K. V. Namboothiri, Certain weighted averages of generalized Ramanujan sums, {\it Ramanujan J.}  {\bf 44} (2017), 531--547.

\bibitem{Ni} N. Nielsen, {\it Handbuch der Theorie der Gammafunktion}, Chelsea Publishing Company. Bronx, New York (1965). 


\bibitem{R.R1} N. Robles and A. Roy, Moments of averages of generalized Ramanujan sums, {\it Monatsh Math.} {\bf 182} (2017), 433-461.

\bibitem{Sin} J. Singh, Defining power sums of $n$ and $\phi(n)$ integers, {\it Int. J. Number Theory} {\bf 5} (2009), 41-53. 

\bibitem{To1} L. T\'{o}th, Some Remarks on a Paper of V. A. Liskovets,
{\it Integers} {\bf 11: A51} (2011).

\bibitem{To2} L. T\'{o}th, Averages of Ramanujan sums: Note on two papers by E. Alkan, {\it Ramanujan J.} {\bf 35} (2014), 149-156.

\bibitem{To3} L. T\'{o}th, Weighted gcd-sum functions, 
{\it J. Integer Sequences} {\bf 14} (2011), Article 11.7.7. 


\end{thebibliography}
\end{document}